\documentclass{article}
\usepackage{graphicx}
\usepackage{amsmath,amssymb,amsfonts,amsthm,amsbsy}
\usepackage{hyperref}
\usepackage{enumerate}

\newtheorem{definition}{Definition}[section]
\newtheorem{proposition}[definition]{Proposition}
\newtheorem{theorem}[definition]{Theorem}
\newtheorem{corollary}[definition]{Corollary}
\newtheorem{example}[definition]{Example}
\newcommand{\overbar}[1]{\mkern 1.5mu\overline{\mkern-1.5mu#1\mkern-1.5mu}\mkern 1.5mu}
\def\nnabla{\nabla \hskip-2.7mm \nabla}

\begin{document}
\title{Symplectic connections and Fedosov's quantization on supermanifolds}
\author{Jos\'e A. Vallejo \\
{\normalsize Facultad de Ciencias, Universidad Aut\'onoma de San Luis Potos\'i (M\'exico)}\\
{\footnotesize Email: \texttt{jvallejo@fc.uaslp.mx}}
}
\date{\today}
\maketitle

\begin{abstract}
A (biased and incomplete) review of the status of the theory of symplectic connections on supermanifolds is presented. Also, some comments regarding Fedosov's technique of quantization are made.
\end{abstract}

\section{Introduction}
The quantization problem can be described as the search for a mapping from the classical observables of a physical system
(smooth functions on a symplectic or, more generally, Poisson manifold) to self-adjoint operators on a certain Hilbert space $\mathcal{H}$ (a space which must be constructed during the process). In the early days of Quantum Mechanics, the naive procedure due to Dirac (\cite{Dir-58}) was followed, basically consisting in the replacements for coordinates and momenta $q \mapsto \hat{q}\psi =x\cdot \psi$ (the hat denoting operators, and $\psi \in \mathcal{H}$) and $p\mapsto \hat{p}\psi =-i\frac{\mathrm{d}\psi}{\mathrm{d}q}$.\\
This simple scheme does not work for arbitrary systems (notably for those classical Hamiltonians containing cubic or higher order mixed terms in $p,q$, see \cite{Gro-46,VH-51,Cher-81}), due to the fact that there is not a unique way to assign a product of operators to a product of commutative expressions such as $p^2 q,pq^2$, etc. Weyl tried to solve this problem within the ``phase-space formalism'' of Quantum Mechanics (see \cite{ZFC-05} and references therein for an overview of this topic) by fixing an ordering, establishing a correspondence of the form
$$
q^m p^n \mapsto \frac{1}{2}\sum^{n}_{k=0} {a \choose b} \hat{p}^{n-k} \hat{q}^m \hat{p}^k .
$$
The next step was the generalization of this construction to classical functions other than polynomials. This was accomplished by Weyl through the use of a certain class of pseudo-differential operators (see \cite{Wey-50,Shu-01}): if we have a classical observable $f(q,p)$ in $\mathbb{R}^{2n}$, which is to be thought as the phase space of the system, it can be formally expressed, by using the inverse Fourier transform, as
$$
\hat{f}(\xi,\eta)= \iint_{\mathbb{R}^{2n}} e^{-\frac{i}{\hbar}(q\xi +p\eta)}f(q ,p)\mathrm{d}q \mathrm{d}p ,
$$
and
\begin{equation}\label{fourier-trans}
f(q,p)= (2\pi \hbar)^{-2n}\iint_{\mathbb{R}^{2n}}e^{\frac{i}{\hbar}(q\xi +p\eta)}\hat{f}(\xi ,\eta )\mathrm{d}\xi \mathrm{d}\eta .
\end{equation}
Weyl then defines the operator corresponding to $f$ as the $W_f$ obtained by substituting $q,p$ for $\hat{q},\hat{p}$ in the formula \eqref{fourier-trans}, so
$$
W_f = (2\pi \hbar)^{-2n}\iint_{\mathbb{R}^{2n}} e^{\frac{i}{\hbar}(\xi \hat{q} +\eta \hat{p})}\hat{f}(\xi ,\eta )\mathrm{d}\xi \mathrm{d}\eta ,
$$
and this acts on a function $u\in \mathcal{L}^2 (\mathbb{R}^n )$ through its kernel, giving\footnote{After some algebraic manipulations, the use in the Baker-Campbell-Hausdorff formula of the commutator $\left[ p\frac{\mathrm{d}}{\mathrm{d}q},i\eta q\cdot \right] = i\eta p\mathbb{I}$, and the fact that $e^{i\partial_q}$ is the generator of translations in $q$.} the result:
$$
(W_f u)(q)=(2\pi \hbar )^{-2n}\iint_{\mathbb{R}^{2n}}e^{\frac{i}{\hbar}(q-p)\xi}f\left( \frac{q+p}{2},\xi \right)u(p) \mathrm{d}\xi \mathrm{d}p .
$$
Let us remark that $W_f$ acts on functions defined on $\mathbb{R}^n$ (the configuration space) and not on the full phase-space $\mathbb{R}^{2n}$. The mathematical details justifying these formal manipulations can be found in \cite{Gad-95,GBV-98}. What we want to remark now, is the fact that the correspondence $f\mapsto W_f$ is $\mathbb{C}-$linear, but the space of classical observables is not just a vector space, it also has the structure of a commutative algebra, with the point-wise product of functions. The space of self-adjoint operators such as $W_f$, when endowed with the composition of operators, also has an algebra structure, but this time a non-commutative one. To describe Quantum Mechanics in phase space following the initial motivation, one should be able to establish also a morphism between these two algebras, but this can not be done without modifying the commutative product of classical observables. Thus, the problem was to find a non-commutative product $\star$, on the algebra of functions on $\mathbb{R}^{2n}$, such that the operator corresponding to $f\star g$ is precisely $W_f \circ W_g$. Clearly, this can be done by defining
\begin{equation}\label{weyl-prod}
f\star g=W^{-1}(W_f \circ W_g).
\end{equation}
The expression for the inverse operator $W^{-1}$ was found by Wigner short after the ideas of Weyl were published, in 1932 (see \cite{Wig-32}).
However, the formulae found by Wigner (in terms of a trace operator) did not allow for a direct physical interpretation, so J. E. Moyal undertook the task of finding that interpretation for $f\star g$. He (and, independently, H. Groenewold) discovered what nowadays we call the Moyal product (see \cite{Gro-46,Moy-49}). In modern terminology, Moyal realized that the star product can be written in terms of the Poisson bi-vector determined by the classical bracket on $\mathcal{C}^{\infty}(\mathbb{R}^{2n})$,
$$
\{ f,g \}=P(f\otimes g)=\sum \left( \frac{\partial f}{\partial q^i}\frac{\partial g}{\partial p_j}-\frac{\partial f}{\partial p_j}\frac{\partial g}{\partial q^i}\right),
$$
as
\begin{equation}\label{moyal-prod}
f\star_{\hbar} g=\mathrm{exp}\left( \frac{i\hbar}{2}P\right) (f\otimes g)=f\cdot g+\sum^{\infty}_{k=1}\left( \frac{i\hbar}{2} \right)^k P^k (f\otimes g).
\end{equation}
Here, $P^k$ denotes the $k-$th iteration of the bi-differential operator $P=\sum \partial_q \wedge \partial_p$. Thus, the classical Poisson bracket can be recovered as the limit $\hbar \to 0$ of the so-called Moyal bracket:
$$
[ f,g]_{\hbar} =\frac{1}{i\hbar}(f\star_{\hbar} g-g\star_{\hbar} f).
$$
Indeed:
$$
\{ f,g \}=\lim_{\hbar \to 0}[f,g]_{\hbar} .
$$
In this sense, the Moyal bracket plays the role of the quantum bracket (but using only quantities defined in the classical phase space, instead of operators, hence the name ``Quantum Mechanics on phase space''), and the classical Poisson bracket is the limit when $\hbar \to 0$ of the quantum bracket.\\
This idea can be reversed, and then one can study star products as \emph{deformations} of the point-wise product of observables defined on classical phase space, trying to get in this way some information about the structure of quantum mechanical systems, as originally proposed by Dirac in \cite{Dir-58}. This was also the proposal by Bayen, Flato, Fr\o nsdal, Lichnerowicz and Sternheimer in a series of papers that laid the foundations of deformation quantization (see \cite{Fla-74,FLS-74,FLS-76,BFFLS-77,BFFLS-78}). We will not deepen here into the details of this interesting formalism, which has attracted a lot of attention, as there exist several excellent reviews in the literature (see for example \cite{Wei-94,DS-02,Bor-08} and references therein).

Let us comment the contents of the paper. This introduction intends to put the study of Fedosov formalism in its historical and physical context.
Section \ref{sec-2} describes the basics of Fedosov's construction putting the emphasis on the r\^{o}le played by the symplectic connections, which are considered in greater detail in Section \ref{sec-3}. The geometry of supermanifolds, including connections and symplectic forms, is recalled in Section \ref{sec-4} and then, in the last Section, we briefly report on the various approaches that have been developed to study Fedosov supermanifolds, comparing the results obtained and remarking the remaining open questions. Although this survey is by no means complete, an effort has been done to provide precise references studying in detail the topics which are just touched upon here, while some attention has been paid to technical details often glossed over.

\section{Fedosov's quantization}\label{sec-2}
The Moyal-Weyl theory is very nice, but has a serious problem. Usually, the phase spaces of dynamical systems have the structure of (curved) manifolds, not simply Euclidean spaces. However, it is not evident at all how to generalize the Moyal-Weyl formalism to this case. Of course, every manifold can be constructed by gluing together Euclidean domains, but if one wants to obtain a global Moyal product by patching local ones, the precise procedure becomes very intricate. It was not until recently that this construction was done for the case of symplectic manifolds, independently by De Wilde and Lecomte \cite{WL-83}, Omori, Maeda and Yoshioka \cite{OMY-91}, and Fedosov \cite{Fed-94}. For Poisson manifolds, using different techniques based on the study of the Formality Conjecture \cite{Kon-97a}, the goal was achieved by Kontsevich \cite{Kon-97}. Cattaneo and Felder (in \cite{CF-01}) gave a connection between the methods of Kontsevich and those of string theories (the Poisson-sigma model), globalizing some local constructions in \cite{Kon-97} \`a la Fedosov. Here we will be interested in Fedosov's approach, which has a strong geometrical flavour and is the most used in Physics.

Thus, the goal of deformation quantization is to provide a star product with which a Moyal-type bracket can be defined, having as the formal limit when $\hbar \to 0$ the classical Poisson bracket. We now give the precise definition.
\begin{definition}\label{def-1}
Let $(M,\{ \cdot ,\cdot \})$ be a Poisson manifold. A (differential) star product $\star$ on $(M,\{ \cdot ,\cdot \})$ is a formal associative $\mathbb{C}[[\hbar ]]-$bilinear product defined on the space $\mathcal{C}^{\infty}(M)[[\hbar ]]$, given by
$$
f\star g =\sum^{\infty}_{j=0}C_j (f,g)\hbar^j ,
$$
where the $C_j$ are bi-differential operators on $\mathcal{C}^{\infty}(M)[[\hbar ]]$ such that
\begin{enumerate}[(a)]
\item $C_0 (f,g)=f\cdot g$ (the point-wise product).
\item $C_1 (f,g)-C_1 (g,f)=i\hbar \{ f,g \}$.
\item $\mathbf{1}\star f =f\star \mathbf{1}$.
\end{enumerate}
A (formal) deformation quantization on $(M,\{ \cdot ,\cdot \})$ is a pair $(\mathcal{C}^{\infty}(M)[[\hbar ]],\star )$. The quantization is said to be strict when $\hbar$ is not just a formal parameter, but a real one.
\end{definition}

Fedosov's constructs $\mathcal{C}^{\infty}(M)[[\hbar ]]$ as the center of the algebra of sections of a certain bundle, called the Weyl bundle. Some of the technical notions and results involved in his construction will be discussed in the next section, but the idea is the following. Take a manifold $N$ and a coordinate chart $( U,\{ q^i \})$ on $N$, so we have induced basis $( TU,\{ q^i, \frac{\partial}{\partial q^i} \})$ on $TN$ and  $( T^{\ast}U,\{ q^i, \mathrm{d}q^i \})$ on $T^{\ast}N$. Let $\hbar$ be a positive parameter, and consider at each point $p\in T^{\ast}U$ the space of formal series of the form
$$
a =a(\hbar, X)=\sum^{\infty}_{l=0} \hbar^k \alpha_{k,i_1 ,...,i_l}X^{i_1}\cdots X^{i_l}, \quad k\geq 0,
$$
where $X^{i_1},...,X^{i_l}$ are the components of a vector $X\in T_p (T^{\ast}N)$ and $\alpha_{k,i_1 ,...,i_l}$ are the components of a covariant tensor on $M$ symmetric with respect to indices $(i_1 ,...,i_l )$ (in the basis $\mathrm{d}q^{i_1} \otimes \cdots \otimes \mathrm{d}q^{i_l}$). If we write $\mathcal{V}_p$ for the set of all the formal series $a$ at the point $p\in  T^{\ast}U$, it has the structure of a vector space, indeed, a Fedosov's product $\circ$ can be defined on $V_p$ so it becomes a non-commutative algebra. This product is an extension of the Moyal-Weyl product on $\mathbb{R}^{2n}$ and the exterior product:
\begin{equation}\label{fed-prod}
a\circ b = \sum^{\infty}_{r=0}\left( \frac{i\hbar}{2} \right)^r \frac{1}{r!}\omega^{i_1 j_1}\cdots \omega^{i_r j_r}
\frac{\partial^r a}{\partial X^{i_1}\cdots \partial X^{i_r}} \cdot \frac{\partial^r b}{\partial X^{j_1}\cdots \partial X^{j_r}},
\end{equation}
where $\omega^{ij}$ are the components of the inverse of the canonical symplectic form $\omega$ on $T^{\ast}N$. Then we can form the bundle
\begin{equation}\label{bundle-v}
\mathcal{V}=\bigcup_{p\in T^{\ast}N}\mathcal{V}_p ,
\end{equation}
which is the bundle of forms valued in the Weyl bundle of $TN$.\\
This construction can be seen independent of the chosen charts when written in geometric terms. To this end, let us directly consider a symplectic manifold $(M,\omega )$ \emph{endowed with a symplectic connection} $\nabla$, that is, such that $\nabla \omega =0$. On each  symplectic tangent space $(T_p M,\omega_p )$ (which is isomorphic to $\mathbb{R}^{2n}$) we can define the Moyal-Weyl star product $\star_{\hbar}$ as follows. There is a unique differential operator $P:\mathcal{C}^{\infty}(T_p M\times T_p M)\to \mathcal{C}^{\infty}(T_p M\times T_p M)$, with constant coefficients, such that
\begin{equation}\label{eq-delta}
\{ f,g \} =\Delta^{\ast} P(f\otimes g),
\end{equation}
where $\{ \cdot ,\cdot \}$ is the (non-degenerate) Poisson bracket induced by $\omega$, $f\otimes g$ is given by $f\otimes g (u,v)=f(u)g(v)$ for any $u,v\in T_p M$, and $\Delta^{\ast}:\mathcal{C}^{\infty}(T_p M\times T_p M)\to \mathcal{C}^{\infty}(T_p M)$ is the restriction to the diagonal. The Moyal-Weyl product is then
$$
f\star_{\hbar} g=\Delta^{\ast}\mathrm{exp}\left( \frac{i\hbar}{2} P\right)(f\otimes g).
$$
The space $(\mathcal{C}^{\infty}(T_p M)[[\hbar ]],\star_{\hbar})$ is called the Weyl algebra of $T_p M$, usually denoted by $W(T_p M)$ or simply $W_p (M)$. These algebras determine a fibre bundle over $M$, $W(M)=\cup_{p\in M}W_p (M)$, which is the Weyl bundle of $M$. The space of sections $\Gamma W(M)$ has a natural structure of unital associative algebra induced by the fibre-wise Moyal product, and the center of $\Gamma W(M)$ can be identified (as a $\mathbb{C}-$vector space) with $\mathcal{C}^{\infty}(M)[[\hbar ]]$, as mentioned.\\
Consider now the bundle $\mathcal{V}=\Gamma W(M)\otimes \Omega (M)$ \eqref{bundle-v}, of $\Gamma W(M)-$valued differential forms on $M$, which is endowed with the natural product $\circ$ described before. It is also endowed with the extension to $\Gamma W(M)\otimes \Omega (M)$ of the projection onto the center $\tau : \Gamma W(M)\to \mathcal{C}^{\infty}(M)[[\hbar ]]$, an extension which we will also denote by $\tau$.\\
The symplectic connection on $M$, $\nabla$, induces a connection on $W(M)$, $D_0 :\Gamma W(M)\otimes \Omega^k (M)\to \Gamma W(M)\otimes \Omega^{k+1}$ through the local expression
$$
D_0 =\mathrm{d}x^i \wedge \nabla_i.
$$
Then, for an $a\in \Gamma W(M)\otimes \Omega^1 (M)$ such that $\tau (a)=0$, Fedosov defines the connection
$$
D_1=D_0 -\frac{1}{i\hbar}[\![ a,\cdot ]\!],
$$
where $[\![ \cdot ,\cdot ]\!]$ is the graded bracket $[\![ a,b]\!] =a\circ b-(-1)^{|a||b|}b\circ a$ for $a\in \Gamma W(M)\otimes \Omega^{|a|} (M)$ and $b\in \Gamma W(M)\otimes \Omega^{|b|} (M)$. This new connection is not flat, in the sense that $D^2_1 \neq 0$, so Fedosov gives an iterative method to modify it in order to get a new connection $D$ (called the Fedosov connection) such that $D^2 =0$.\\
Once we have the flat connection $D:\Gamma W(M)\otimes \Omega^k (M)\to \Gamma W(M)\otimes \Omega^{k+1}$ (which, we insist, is algorithmically  constructed from $\nabla$), the set of its horizontal sections, that is, those $\sigma$ for which $D\sigma =0$, forms an algebra with respect to the product $\circ$, a subalgebra that will be denoted $(\Gamma_D (W),\circ )$ or simply $\Gamma_D (W)$. The restriction of the projection map $\tau :\Gamma W(M)\to \mathcal{C}^{\infty}(M)[[\hbar ]]$ to $\Gamma_D (W)$ is bijective so, to any $f\in \mathcal{C}^{\infty}(M)[[\hbar ]]$, we can associate a horizontal section $\sigma_f =\tau^{-1}(f)$. Then,
\begin{equation}\label{fedosov-star}
f\star_{\hbar} g =\tau (\sigma_f ,\sigma_g )
\end{equation}
defines a star product on $M$ and hence a deformation quantization of $(M,\omega )$, in the sense of Definition \ref{def-1}. Note the similarity between \eqref{fedosov-star} and \eqref{weyl-prod}. Indeed, for $M=\mathbb{R}^{2n}$ with Euclidean coordinates, the product \eqref{fedosov-star} is just the Moyal product \eqref{moyal-prod}, where $P$ is the (non-degenerate) Poisson bi-vector determined by the symplectic form $\omega$. 

\section{Symplectic connections and Fedosov manifolds}\label{sec-3}
In the preceding section, we have tried to stress the r\^ole of the symplectic connections $\nabla$ in Fedosov's approach, now we would like to explore some of the properties of these objects. A symplectic connection is just a torsion-free linear connection on the symplectic manifold $(M,\omega )$, such that the symplectic form is parallel. Let us recall the basic formalism. 
\begin{definition}\label{defconnection}
Let $E$ be a vector bundle over a manifold $M$. A linear connection $\nabla$ on $E$ is an $\mathbb{R}-$linear mapping $\nabla :\Gamma E \to \Omega^1 (M;E)$ (vector-valued differential forms on $M$) such that
$$
[\nabla ,\mu_f ]=\mathrm{d}f,
$$
where, given a $f\in \mathcal{C}^{\infty}(M)$, $\mu_f$ denotes the product by $f$ and $[ \cdot ,\cdot ]$ is the commutator of bundle endomorphisms.
\end{definition}
Denoting the action of the connection on a pair $(X,\sigma )\in \Gamma (TM\otimes E)$ by $\nabla_X \sigma$, the torsion $T^{\nabla}$ and curvature $R^{\nabla}$ of $\nabla$ are defined as usual; for the curvature we have:
$$
R^{\nabla}(X,Y)=[\nabla_X ,\nabla_Y ]-\nabla_{[X,Y]},
$$
and, taking $E=TM$, the torsion is given by: 
$$
T^{\nabla}(X,Y)=\nabla_X Y-\nabla_Y X-[X,Y].
$$
When $T^{\nabla}=0$ the connection is said to be torsion-free (or symmetric), and it is flat if $R^{\nabla}=0$. Note that it is equivalent to define a connection on a vector bundle as in Definition \ref{defconnection} or as a derivation of degree $1$ on $\Omega (M,E)=\Gamma (\bigwedge T^{\ast}M\otimes E)$, $D:\Omega^k (M,E)\to \Omega^{k+1} (M,E)$  such that the induced derivation $\overbar{D}:\Omega^k (M)\to \Omega^{k+1} (M)$ given by 
\begin{equation}\label{eq-D}
[D, \alpha ]=\overbar{D}\alpha ,
\end{equation}
is precisely the exterior derivative $\mathrm{d}$. In this case it is common to write $D=d^{\nabla}$ and to call it the covariant exterior derivative (this is precisely the form in which the connections $D_0,D_1,...$ in the Fedosov's construction are given). In fact, given a connection $\nabla$ on $E$, the associated covariant derivative $d^{\nabla}$ is defined (being a derivation of degree $1$) through its action on $\mathcal{C}^{\infty}(M)$ and on $\Omega^0 (M;E)\simeq \Gamma E$:
$$
d^{\nabla}f=\mathrm{d}f,\quad d^{\nabla}\sigma =\nabla \sigma,
$$  
and extending as a derivation to $\Omega (M;E)$. It is readily seen that the curvature of $\nabla$ is related to $d^{\nabla}$ by
$$
\left( d^{\nabla}\right)^2 \gamma =R^{\nabla}\wedge \gamma ,
$$
for any $\gamma \in \Omega (M;E)$, so $\nabla$ is flat if and only if $\left( d^{\nabla}\right)^2 =0$.\\
From \eqref{eq-D} and the fact that the algebra $\Omega (M;\mathrm{End}(E))$ can be identified with the algebra of operators on
$\Omega (M;E)$ that are $\Omega (M)-$linear, we have:
\begin{proposition}\label{affine}
The space of connections on a vector bundle $E$ over $M$, is an affine space modelled on the vector space $\Omega^1 (M;\mathrm{End}(E))$.
\end{proposition}
The main definition is then:
\begin{definition}\label{def-fedosov}
A symplectic connection on a symplectic manifold $(M,\omega )$ is a linear connection $\nabla$ on $TM$ such that it is torsion-free and $\nabla \omega =0$. A Fedosov manifold is a triple $(M,\omega ,\nabla )$ with $\nabla$ symplectic.
\end{definition}
Up to here, the symplectic and Riemannian connections share the same properties. Indeed, the existence of symplectic connections on any symplectic manifolds is proved exactly in the same way as in the Riemannian case (define them locally in Darboux charts and then glue using a partition of unity). However, there are important differences. For instance, while once given a Riemannian metric $g$ on $M$ there is only one Levi-Civit\`a connection (torsion-free and such that $\nabla g=0$), there are infinite symplectic connections for a given symplectic form $\omega$ and, indeed, they form an affine space (compare with the generic case described in Proposition \ref{affine}):
\begin{proposition}{\cite{BC-99,BCGRS-05}}
Let $(M,\omega )$ be a symplectic manifold. The set of symplectic connections
on $(M,\omega )$ can be identified with the affine space of completely symmetric tensor fields of
$3-$covariant type on $M$.
\end{proposition}
\noindent The proof consists in noticing that given $\nabla$ symplectic, then $\nabla '_X Y=\nabla_X Y+A(X,Y)$ is symplectic if and only if $\omega (A(X,Y),Z)$ is symmetric in all its indices.

That this kind of connections is important in deformation quantization had been already pointed; in fact, in the seminal paper \cite{BFFLS-77} the authors observe that a Moyal star product could be defined on any symplectic manifold $(M,\omega )$ which admits a flat symplectic connection $\nabla$, and A. Lichnerowicz in \cite{Lic-82} remarked that for certain star products $\star$ on a symplectic manifold, there exist a unique symplectic connection satisfying some set of conditions. This has been generalized by S. Gutt and J. Rawnsley in \cite{GR-03}, where they prove that any natural star product on a symplectic manifold determines a unique symplectic connection. Fedosov's results can be seen as a converse of these, as he provides an explicit construction of a star product starting from a symplectic connection. The applications in Physics are not restricted to deformation quantization, for different examples see \cite{DLS-01,Kry-07,Tos-07,Vas-09} and references therein.

Other interesting properties of symplectic connections refer to their curvature. The symplectic curvature was first systematically studied by Vaisman in \cite{Vai-85} (where most of the formulae below are proved), and later more generally in the context of the geometry of Fedosov manifolds by Gelfand, Retakh and Shubin in \cite{GRS-98}. For more recent works on Fedosov manifolds, see \cite{BC-99,BCGRS-05,Dub-05,Gut-06,EM-08}.

Recall that, in Riemannian geometry, given a metric $g$ and its Levi-Civit\`a connection\footnote{Actually, some of the formulae make sense for arbitrary linear connections on $TM$.} $\nabla$ (which is a linear connection on $TM$), the Riemann curvature tensor $R$ is defined as the following tensor of type $(3,1)$:
$$
R^{\nabla}(X,Y)Z=\nabla_X \nabla_Y Z-\nabla_Y \nabla_X Z -\nabla_{[X,Y]}Z.
$$
A direct consequence of the definition, the Jacobi identity for the Lie bracket of vector fields, and the symmetry of the Levi-Civit\`a connection, is the Bianchi's Identity:
\begin{equation}\label{eq-bianchi}
R^{\nabla}(X,Y)Z+R^{\nabla}(Y,Z)X+R^{\nabla}(Z,X)Y=0.
\end{equation}
Other symmetries of the Riemann tensor are better written in terms of its completely covariant form 
\begin{equation}\label{eq-riemann}
R^g (X,Y,Z,T)=g(R^{\nabla}(X,Y)Z,T).
\end{equation}
These are:
\begin{eqnarray}\label{eq-symmetries}
R^g (X,Y,Z,T)=-R^g (Y,X,Z,T) \nonumber \\
R^g (X,Y,Z,T)=-R^g (X,Y,T,Z) \\
R^g (X,Y,Z,T)=R^g (Z,T,X,Y) \nonumber .
\end{eqnarray}
From the Riemann curvature, the Ricci tensor $\rho$ (or Ricci curvature) is constructed as a trace with respect to $g$; it is the $(2,0)$ tensor given by:
$$
\rho (X,Y)=\mathrm{Tr}_g (Z \mapsto R^{\nabla}(X,Z)Y ).
$$
Due to the Bianchi's Identity \eqref{eq-bianchi}, the Ricci tensor is symmetric: $\rho (X,Y)=\rho (Y,X)$.\\
Let us remark that there is another possibility for forming a $2-$covariant tensor from the Riemann curvature $R$. It consists in taking
$$
\tilde{\rho} (X,Y)=\mathrm{Tr}(Z \mapsto R^{\nabla}(X,Y)Z ),
$$
but, from the properties \eqref{eq-bianchi} and \eqref{eq-symmetries}, it is readily seen that this leads to $\tilde{\rho}=0$.\\
By a further contraction of the Ricci tensor, the (Riemannian) scalar curvature $S$ is defined:
$$
S=\mathrm{Tr}_g (Z \mapsto \rho(X,Z)).
$$
Now, given a Fedosov manifold $(M,\omega ,\nabla )$, let us define the analogue of \eqref{eq-riemann}.
\begin{definition}
The symplectic curvature of $(M,\omega ,\nabla )$ is the $(4,0)-$tensor given by
$$
R(X,Y,Z,T)=\omega (R^{\nabla}(X,Y)Z,T).
$$
\end{definition}
Of course, the Bianchi's Identity \eqref{eq-bianchi} is true in this setting (as it only depends on the symmetry of the connection). But we would like to stress the differences between the Riemannian and symplectic curvatures. While the former is skew-symmetric under the transformation \eqref{eq-symmetries} of the indexes, we now have \cite{Vai-85,GRS-98}:
\begin{proposition}\label{pro-propcurv}
The symplectic curvature has the symmetries
\begin{enumerate}[(1)]
\item $R(X,Y,Z,T)=R(Y,X,Z,T)$, for all $X,Y,Z,T\in\mathcal{X}(M)$,
\item $R(X,Y,Z,T)=-R(X,Y,T,Z)$,
\item $R(X,Y,Z,T)+R(T,X,Y,Z)+R(Z,T,X,Y)+R(Y,Z,T,X)=0$.
\end{enumerate}
\end{proposition}

If $(M,\omega )$ is a symplectic manifold with $\mathrm{dim}M=2n$, the non-degenerate symplectic tensor $\omega$ induces a bundle isomorphism (the Poisson morphism) $\mathbf{P}:TM\to T^{\ast}M$. Then, we can define for $\alpha ,\beta \in \Gamma T^{\ast}M$:
$$
\omega^{-1} (\alpha ,\beta )= \omega (\mathbf{P}^{-1} \alpha ,\mathbf{P}^{-1} \beta ).
$$
If a set of vector fields $\{ U_i\}^{2n}_{i=1}$ is a local frame, we will consider the dual frame $\{ \alpha^j \}^{2n}_{j=1}$, such that
$\alpha^j (U_i )=\delta^j_i$.

With respect to the Ricci tensor, we have the following
\begin{theorem}[\cite{GRS-98}]\label{teo-contractions}
Suppose $(M,\omega ,\nabla )$ is a Fedosov manifold with $\mathrm{dim}M=2n$. There are two possible non-trivial contractions of the symplectic curvature:
$$
\rho (X,Y)=\mathrm{Tr}_{\omega}(Z\mapsto R^{\nabla}(Z,Y)X),
$$
\begin{equation}\label{rotilde}
\tilde{\rho}(X,Y)=\sum_{i,j}\omega^{-1} (\alpha^j ,\alpha^k )\omega (R^{\nabla}(U_j ,U_k )X,Y). 
\end{equation}
The relationship between them is
$$
\tilde{\rho}=2\rho .
$$
Moreover, $\rho$ (and so $\tilde{\rho}$) is symmetric.
\end{theorem}
\noindent Notice that the contraction given in \eqref{rotilde} is trivial in the Riemannian case, due to the fact that $g^{-1}$ is symmetric while
$R^{\nabla}$ is skew-symmetric, so it gives identically zero.
\begin{corollary}\label{cor-szero}
The symplectic scalar curvature on a Fedosov manifold (the contraction of the Ricci tensor with the symplectic form) vanishes identically.
\end{corollary}
As we will see, this is not the case for Fedosov supermanifolds, where an odd symplectic scalar curvature appears.

\section{Geometry on supermanifolds}\label{sec-4}
There are several approaches to the notion of supermanifold. Roughly speaking, they can be classified in differential-geometric (\cite{Wit-92,Rog-86,Rog-07,JP-81}) and sheaf-theoretic (\cite{Kos-77,Lei-80,Man-97}) ones (for expositions of the theory oriented towards physics applications, see \cite{Var-04,Sar-09}). Also, there are formulations of the theory of supermanifolds that takes advantages of the best of both approaches (\cite{BBH-91,Tuy-04}). In some cases, these are equivalent (see \cite{Rot-86,BBH-91}) or at least there is a bijective correspondence (not functorial) between different categories of supermanifolds (\cite{BP-99}), but there are technical details that must be taken into account when choosing one or another.

Let us describe first the differential-geometric approach. In it, one begins by fixing a Grassmann algebra $B$, which splits into its even and odd parts, $B=B_0 \oplus B_1$. Consider now the following space $B^{m,n}$, parametrized by $m$ commuting and $n$ anticommuting variables:
$$
B^{m,n}=B_0 \times \overset{m)}{\cdots} \times B_0 \times B_1 \times \overset{n)}{\cdots} \times B_1 .
$$
It is assumed that this space $B^{m,n}$ is endowed with some topology; for instance, in \cite{Rog-80} this is done seeing it as a Banach space, in 
\cite{IM-91} as a Fr\'echet space, and in \cite{Wit-92} through a non-Hausdorff topology. Now, to define a geometric supermanifold of dimension $(m,n)$, the usual definition of charts and atlases is used, but replacing the Euclidean model $\mathbb{R}^n$ as the target for the chart mappings, by the local model $B^{m,n}$; thus, a geometric $(m,n)-$dimensional supermanifold is a topological space $M$ covered b y an atlas $\{ V_{\alpha},\varphi_{\alpha}\}_{\alpha \in I}$ such that $\varphi_{\alpha}:V_{\alpha}\to B^{m,n}$. The main point here is that the transition functions must be superdifferentiable in some appropriate sense, which depends on the nature of the Grassmann algebra $B$. In the original formulation of A. Rogers \cite{Rog-80}, $B$ is finite-dimensional, and the definition of superdifferentiability for functions $f:B^{m,n}\to B^{m,n}$ (the so-called $G^{\infty}-$differentiability) is flawed because there is not a unique operator for partial differentiation along odd directions \cite{BG-84}. Indeed, M. Rothstein proved in \cite{Rot-86} that, in order to have a structure of a locally free module on the module of derivations of $G^{\infty}-$superfunctions, the number of odd generators of $B$ must be $0$ or infinite. The study of $G^{\infty}-$supermanifolds based on an infinite number of odd generators is done in \cite{JP-81}. In a later work \cite{Rog-86}, A. Rogers introduced the $GH^{\infty}-$differentiability, constructing $GH^{\infty}-$supermanifolds with an arbitrary (finite or infinite) number of odd generators. But this new category has also problems, because it is not possible to construct in it a sensible notion of supervector bundles in such a way that they are finitely generated and locally free, see \cite{BB-87}. On the other hand, for infinite-dimensional Banach-Grassmann algebras with trivial annihilator (that is, for each odd $a\in B_1$ there exists an odd $b\in B_1$ such that $ab\neq 0$), it is possible to construct a real differentiable function which does not extend to a superdifferentiable function on $B^{m,n}$, see \cite{Pes-93}.\\
A class of geometric supermanifolds that circumvent these difficulties is that of $G-$supermanifolds \cite{BBH-91}, defined as pairs $(M,\mathcal{A})$, where $M$ is a geometric supermanifold modelled on a finite-dimensional Grassmann algebra $B$ and $\mathcal{A}$ is the sheaf given by the tensor product $\mathcal{SC}^{\infty}_M \otimes B$, where  $\mathcal{SC}^{\infty}_M$ is the sheaf of superdifferentiable functions on $M$. This definition is designed to allow for the construction of supervector bundles having the same algebraic and geometric properties as their classical (non-super) counterparts.\\
Indeed, $G-$supermanifolds represent a compromise between geometric supermanifolds and the sheaf-theoretic approach, as we will presently see. In the latter, a supermanifold (also called a graded manifold) is defined as follows.
\begin{definition}\label{def-smfd}
A real supermanifold is a ringed space $(M,\mathcal{A})$, where $\mathcal{A}$ is a sheaf of $\mathbb{Z}_2 -$graded commutative $\mathbb{R}-$algebras such that:
\begin{enumerate}[(a)]
\item If $\mathcal{N}$ denotes the sheaf of nilpotents of $\mathcal{A}$, then $\mathcal{A}/\mathcal{N}$ induces on $M$ the structure of a differential manifold.
\item\label{itemb} The subsheaf $\mathcal{N}/\mathcal{N}^2$ is a locally free sheaf of modules, with $\mathcal{A}$ \emph{locally} isomorphic to the exterior sheaf $\bigwedge \left( \mathcal{N}/\mathcal{N}^2 \right)$.
\end{enumerate}
\end{definition}
Holomorphic graded manifolds are defined similarly, by using sheaves of $\mathbb{C}-$algebras.
\begin{example}\label{ex-koszul}
It is useful to consider the so-called Koszul supermanifold \cite{Kos-94} $(M,\Omega (M))$ as an example of this approach. Here $\Omega (M)=\bigoplus_{p\in \mathbb{Z}}\Omega^p (M)$ is the Cartan exterior algebra formed by the exterior differential forms of the manifold $M$. The nilpotents are all the $\alpha \in \Omega^p (M)$ with $p\geq 1$, so $\mathcal{A}/\mathcal{N}=\mathcal{C}^{\infty}(M)$ (the smooth functions on $M$) in this case. Moreover, $\mathcal{N}/\mathcal{N}^2=\Omega^1 (M)$, the space of $1-$forms, locally generated by the differentials $\mathrm{d}x^1 ,...,\mathrm{d}x^m$ of the functions $x^i$ of a chart on $M$.
\end{example}
\noindent From this example, it is clear that we can think of a supermanifold $(M,\mathcal{A})$ as a classical manifold $M$ where the sheaf of commutative algebras $\mathcal{C}^{\infty}(M)$ has been replaced by a sheaf of $\mathbb{Z}_2 -$graded commutative $\mathbb{R}-$algebras, $\mathcal{A}$, whose sections are the superfunctions.

The local isomorphism in \eqref{itemb} of Definition \ref{def-smfd}, means that there exists a covering of $M$ by open sets $U\subset M$ such that
\begin{equation}\label{sfunctions}
\mathcal{A}(U)\simeq \mathcal{C}^{\infty}(U)\otimes \bigwedge (\mathbb{R}^n ).
\end{equation}
The dimension of $(M,\mathcal{A})$ as a supermanifold is then $(m,n)$, where $m=\mathrm{dim}M$. The open sets $U$ are called splitting neighbourhoods. In them, a local set of supercoordinates is given by $\{ x^i ,\theta_a \}$ (where $i\in \{ 1,...,m\}$, $a\in \{1,...,n\}$), $\{ \theta_1 ,...,\theta_n \}$ being a set of generators for $\bigwedge (\mathbb{R}^n )$. A fundamental theorem of M. Batchelor \cite{Bat-79} states that for a real graded manifold, letting
$\mathcal{E}=\mathcal{N}/\mathcal{N}^2$, we have 
\begin{equation}\label{splitting}
\mathcal{A}\simeq \bigwedge \mathcal{E}
\end{equation}
\noindent not only locally, but globally. However, this sheaf isomorphism is \emph{not} canonical. When the identification \eqref{splitting} is made, it is said that the supermanifold $(M,\mathcal{A})$ is given in split form. An interesting result, due to Koszul \cite{Kos-94}, is that a graded manifold splits if and only if $\mathbb{Z}_2 -$graded connections exist (see below for the definition and main properties of graded connections). This is always the case for real graded manifolds (by Batchelor's theorem), but not for holomorphic ones, see \cite{Gre-82}. Thus, when dealing with graded connections, as we will do, we can always assume that our supermanifold is split\footnote{Note that this is not a loss of generality: ordinary superspaces, when constructed over non trivial finite ground algebras, are the images of split graded manifolds $(M,\bigwedge \mathcal{E})$ under Weil's functor of points, see \cite{BP-99,Var-04,CCF-11}.}. 

From \eqref{sfunctions}, it is immediate that in a splitting neighbourhood  $U\subset M$, with $\{ x^i \}^m_{i=1}$ some chart coordinates and $\{ \theta_a \}^n_{a=1}$ a system of generators of $\bigwedge (\mathbb{R}^n )$, the space of derivations $\mathrm{Der}_{\mathbb{R}}\mathcal{A}(U)$ has the structure of a finitely generated locally free $\mathcal{A}(U)-$module. This is taken as the space of local supervector fields (the $G$ stands for ``graded'') $\mathcal{X}^G (U)$. A system of generators for these supervector fields is then
$$
\left\lbrace \frac{\partial}{\partial x^i},\frac{\partial}{\partial \theta_a} \right\rbrace^{1\leq a \leq n}_{1\leq i \leq m},
$$ 
where
\begin{center}
\begin{tabular}{cc}
$\displaystyle \frac{\partial}{\partial x^i}(x^j)=\delta^j_i$, & $\displaystyle\frac{\partial}{\partial x^i}(\theta_a)=0$, \\
\\
$\displaystyle\frac{\partial}{\partial \theta_a}(x^i)=0$, & $\displaystyle\frac{\partial}{\partial \theta_b}(\theta_a)=\delta^a_b$. \\ 
\end{tabular}
\end{center}
Of course, the notation $\mathcal{X}^G (M)$ will mean $\mathrm{Der}_{\mathbb{R}}\mathcal{A}$. In the case of geometric supermanifolds, as mentioned before, some problems appear when constructing the supertangent space (and, in general, supervector bundles) in such a way that it is locally free and finitely generated. The only categories of geometric supermanifolds free from these problems are DeWitt and $G-$supermanifolds. DeWitt supermanifolds are equivalent to graded manifolds (see \cite{BBH-91}, Chapter V) and, in turn, $G-$supermanifolds are related to DeWitt ones by an extension procedure. As a consequence, the tangent spaces to these have basically the same structure.

Once we have defined supervector fields, superforms can be constructed in the usual way, by taking duals and exterior products (taking into account the $\mathbb{Z}_2-$grading, of course); the sheaf of exterior superforms will be denoted $\Omega^G (M)$. Other constructions in classical differential geometry are analogously transferred to the super setting. Let us consider the example of symplectic forms, with their associated Poisson brackets, and connections.

\subsection{Symplectic forms}
As we will be interested in supermanifolds admitting a connection, we will assume a split supermanifold $(M,\bigwedge \mathcal{E})$ (but notice that this is just for simplicity: there will be nothing depending on the splitting, and everywhere we can replace $\bigwedge \mathcal{E}$ by $\mathcal{A}$).
\begin{definition}
A graded $p-$form $\boldsymbol{\alpha}$ on $(M,\Gamma (\bigwedge \mathcal{E}))$ is a $\Gamma (\bigwedge \mathcal{E})-$multilinear alternating morphism
$$
\boldsymbol{\alpha} :\mathcal{X}^G (M)\times \overset{p)}{\cdots}\times \mathcal{X}^G (M)\to \Gamma (\bigwedge \mathcal{E}).
$$
\end{definition}
\noindent The graded $p-$forms form a right $\Gamma (\bigwedge \mathcal{E})-$module $\Omega^p_G (M)$. We also have a graded module $\Omega_G (M)=\bigoplus_{p \in \mathbb{Z}}\Omega^p_G (M)$ (with $\Omega^p_G (M)=\{ 0\}$ for $p<0$). A graded $p-$form is said to have $\mathbb{Z}-$degree $k$ if, for $D_1 ,...,D_p \in \mathcal{X}^G (M)$,
$$
|\boldsymbol{\alpha} (D_1 ,...,D_p )|=\sum^p_{i=1}|D_i |+k.
$$
In this case, it is customary to say that $\boldsymbol{\alpha}$ has bidegree $(p,k)$.
\begin{definition}
The graded exterior differential of a graded $p-$form $\boldsymbol{\alpha} \in \Omega^p_G (M)$, is the graded $(p+1)-$form $\mathbf{d}\boldsymbol{\alpha}\in\Omega^{p+1}_G (M)$ defined by
\begin{eqnarray*}
(\mathbf{d}\boldsymbol{\alpha} )(D_1 ,...,D_{p+1}) &=& \sum^{p+1}_{i=1}(-1)^{i-1+d_{i-1}|D_i |}D_i (\boldsymbol{\alpha} (D_1 ,...,\hat{D_i},...,D_{p+1})) \\
									  &+& \sum_{k<l}(-1)^{d_{k,l}}\boldsymbol{\alpha} ([D_k ,D_l ],D_1 ,...,\hat{D_k},...,\hat{D_l},...,D_{p+1}),
\end{eqnarray*}
where $d_i =\sum^i_{r=1}|D_r |$ and $d_{k,l}=|D_k |d_{k-1}+|D_l |d_{l-1}+|D_k ||D_l |+k+l$. A graded $p-$form $\boldsymbol{\alpha}$ is said to be closed if $\mathbf{d}\boldsymbol{\alpha} =0$.
\end{definition}
\noindent Here, and in what follows, $[ \cdot ,\cdot ]$ denotes the graded commutator of graded endomorphisms, $[E,F]=E\circ F-(-1)^{|E||F|}F\circ E$.\\
Note that, if $\boldsymbol{\alpha}$ has $\mathbb{Z}-$degree $k$, $\mathbf{d}\boldsymbol{\alpha}$ has also $\mathbb{Z}-$degree $k$. It is said that $\mathbf{d}$ is a graded operator of bidegree $(1,0)$.
\begin{example}
If $\boldsymbol{\beta} \in \Omega^1_G (M)$, we have a $2-$form $\mathbf{d}\boldsymbol{\beta} \in \Omega^2_G (M)$ whose action is given by
$$
\mathbf{d}\boldsymbol{\beta} (D_1 ,D_2 )=D_1 (\boldsymbol{\beta} (D_2 ))-(-1)^{|D_1||D_2|}D_2 (\boldsymbol{\beta} (D_1))-\boldsymbol{\beta} ([D_1 ,D_2 ]).
$$
\end{example}
\begin{definition}
A graded symplectic form is a closed graded $2-$form $\boldsymbol{\omega} \in \Omega^2_G (M)$ that is non singular, that is, the $\Gamma (\bigwedge \mathcal{E})-$linear map it induces
\begin{center}
\begin{tabular}{rcl}
$\mathcal{X}^G (M)$ & $\rightarrow$ & $\Omega^1_G (M)$ \\
$D$ & $\mapsto$ & $\iota_D \boldsymbol{\omega}$, 
\end{tabular}
\end{center}
is an isomorphism.
\end{definition}

\noindent Let us note that an arbitrary graded symplectic form $\boldsymbol{\omega}$ decomposes into the sum of its homogeneous parts:
$$
\boldsymbol{\omega} =\sum_{i=0}^n \boldsymbol{\omega}_{(i)},
$$
where each $\boldsymbol{\omega}_{(i)}$ is a graded $2-$form of $\mathbb{Z}-$degree $i$.

\subsection{Poisson brackets}
Let $(M,\Gamma (\bigwedge \mathcal{E}))$ be a graded split manifold. A graded Poisson bracket on $\Gamma (\bigwedge \mathcal{E})$, of $\mathbb{Z}-$degree $k$, is a mapping $[\![ \cdot ,\cdot ]\!]:\Gamma (\bigwedge \mathcal{E})\times \Gamma (\bigwedge \mathcal{E}) \to \Gamma (\bigwedge \mathcal{E})$ satisfying the following conditions, for homogeneous $\alpha,\beta,\gamma \in \Gamma (\bigwedge \mathcal{E})$ with respective $\mathbb{Z}-$degrees $|\alpha |,|\beta |,|\gamma |$:
\begin{enumerate}[(1)]
\item $\mathbb{R}-$bilinearity,
\item $| [\![ \alpha ,\beta ]\!]|=|\alpha |+|\beta |+k$,
\item $[\![ \alpha ,\beta ]\!]=-(-1)^{(\alpha +k)(\beta +k)}[\![ \beta ,\alpha ]\!]$,
\item\label{Leibniz} $[\![\alpha ,\beta \wedge \gamma ]\!]=[\![ \alpha ,\beta ]\!]\wedge \gamma +(-1)^{(\alpha +k)\beta}\beta\wedge [\![\alpha ,\gamma ]\!]$,
\item $[\![ \alpha ,[\![\beta ,\gamma ]\!]]\!]=[\![[\![\alpha ,\beta ]\!],\gamma ]\!]+(-1)^{(\alpha +k)(\beta +k)}[\![ \beta ,[\![\alpha ,\gamma ]\!]]\!]$.
\end{enumerate}
Given a graded Poisson bracket and a section $\alpha \in \Gamma (\bigwedge^{|\alpha |} \mathcal{E})$, the endomorphism $D_{\alpha}$ of $\Gamma (\bigwedge \mathcal{E})$ defined by $D_{\alpha}=[\![\alpha ,\cdot ]\!]$ is a graded derivation (that is, a supervector field) of $\mathbb{Z}-$degree $k+|\alpha |$, due to property \eqref{Leibniz} above, called the Hamiltonian supervector field associated to the superfunction $\alpha$. Then, as in the non graded case, every graded symplectic form $\boldsymbol{\omega}$ determines a graded Poisson bracket, with the same $\mathbb{Z}-$degree, by the expression
$$
[\![ \alpha ,\beta ]\!]=\boldsymbol{\omega} (D_{\alpha},D_{\beta}).
$$
\begin{example}\label{ex-koszul}
Starting from a Poisson manifold $(M,\{ \cdot ,\cdot \})$, it is easy to construct a Poisson bracket of $\mathbb{Z}-$degree $-1$
on the supermanifold of the Example \ref{ex-koszul}, where $\Gamma (\bigwedge \mathcal{E})=\Omega (M)$ (see \cite{Kos-85}). Simply define it on generators as
\begin{align*}
[\![ f,g ]\!] &= 0 \\
[\![ f,\mathrm{d}g ]\!] &= \{ f,g \} \\
[\![ \mathrm{d}f,\mathrm{d}g ]\!] &= \mathrm{d} \{ f,g \} .
\end{align*}
\end{example}
\noindent However, finding even Poisson brackets is much more complicated. As proved in \cite{Rot-91}, every non degenerate even Poisson bracket on $\Gamma (\bigwedge \mathcal{E})$ is determined by the following ``Rothstein data'': a non degenerate Poisson bracket $\{ \cdot ,\cdot \}$ on $M$, a non degenerate metric $g$ on the fibre bundle $E$, and a connection $\nabla$ on $E$, compatible with the metric $g$. A particular class of even Poisson brackets $[\![ \cdot ,\cdot ]\!]$ which are extensions of a classical one $\{ \cdot ,\cdot \}$ (in the sense that the zeroth order term in $[\![ f,g ]\!]$ is precisely $\{ f,g\}$) is studied in \cite{MV-02}. 

\subsection{Connections}
In order to avoid unnecessary and lengthy discussions about arbitrary supervector bundles, we will deal with linear connections \`a la Koszul, exclusively (see \cite{MoS-96,MoS-97,MaS-00}).

\begin{definition}\label{def-connection}
A graded connection on a graded manifold $(M,\mathcal{A})$ is a mapping
\begin{center}
\begin{tabular}{rcl}
$\nnabla :\mathcal{X}^G (M)\times \mathcal{X}^G (M)$ & $\rightarrow$ & $\mathcal{X}^G (M)$ \\
$(D_1 ,D_2 )$ & $\mapsto$ & $\nnabla_{D_1} D_2$ 
\end{tabular}
\end{center}
such that, for all $D_1 ,D_2 ,D_3 \in\mathcal{X}^G (M)$ and $\alpha \in \mathcal{A}$,
\begin{enumerate}[(1)]
\item $\nnabla_{D_1} (D_2 +D_3 )=\nnabla_{D_1} D_2 +\nnabla_{D_1} D_3$,
\item $\nnabla_{(D_1 +D_2 )}D_3 =\nnabla_{D_1} D_3 +\nnabla_{D_2} D_3$,
\item $\nnabla_{\alpha D_1} D_2 =(-1)^{|\alpha ||\nnabla |}\alpha \nnabla_{D_1} D_2$,
\item $\nnabla_{D_1}(\alpha D_2) =D_1 (\alpha )D_2 +(-1)^{|\alpha |(|D_1 |+|\nnabla |)}\alpha \nnabla_{D_1} D_2$.  
\end{enumerate}
A graded connection in $(M,\mathcal{A})$ is called $\mathbb{Z}-$homogeneous of degree $|\nnabla |$ if for any pair of homogeneous derivations $D_1 ,D_2 \in \mathcal{X}^G (M)$, $\nnabla_{D_1} D_2 \in\mathcal{X}^G (M)$ is homogeneous and $| \nnabla_{D_1} D_2 |=|D_1 |+|D_2 |+|\nnabla |$. If
$|\nnabla |\equiv 0\mathrm{mod}2$, it is said that $\nnabla$ is even, and odd if $|\nnabla |\equiv 1\mathrm{mod}2$ (in either case, the connection is said to be $\mathbb{Z}_2 -$graded).\\
The torsion $T$ of $\nnabla$ is the mapping $T:\mathcal{X}^G (M)\times \mathcal{X}^G (M)\to \mathcal{X}^G (M)$ given by
$$
T(D_1 ,D_2 )=\nnabla_{D_1} D_2 -(-1)^{|D_1 ||D_2 |}\nnabla_{D_2} D_1 -[D_1 ,D_2 ].
$$  
\end{definition}
\noindent The graded curvature of a connection is defined as in the non-graded case.
\begin{definition}
Let $\nnabla$ be a graded connection on $(M,\mathcal{A})$. Its graded curvature $R^{\nnabla}$ is given by the mapping
$$
R^{\nnabla}(D_1 ,D_2 )D_3 =[\nnabla_{D_1}, \nnabla_{D_1}]D_3 -\nnabla_{[D_1 ,D_2]}D_3 ,
$$
for all $D_1 ,D_2 ,D_3 \in \mathcal{X}^G (M)$ (so it is a graded tensor of type $(3,1)$).
\end{definition}
\noindent From now on, unless otherwise explicitly stated, we will consider only even connections.
\begin{proposition}[Bianchi's Identity]
Let $\nnabla$ be a torsion-free graded connection and $R^{\nnabla}$ its graded curvature. Then,
\begin{eqnarray*}
R^{\nnabla}(D_1 ,D_2)D_3 &+&(-1)^{|D_1|(|D_2|+|D_3|)}R^{\nnabla}(D_2 ,D_3)D_1 \\
&+&(-1)^{|D_3|(|D_1|+|D_2|)}R^{\nnabla}(D_3 ,D_1)D_2 =0.
\end{eqnarray*}
\end{proposition}
The proof of this basic result is done as in the non graded case, using the fact that for a torsion-free connection
we have:
$$
\nnabla_{D_1}[D_2 ,D_3]-(-1)^{|D_1|(|D_2|+|D_3|)}\nnabla_{[D_2 ,D_3 ]}D_1 =[D_1 ,[D_2 ,D_3 ]],
$$
and applying the graded Jacobi's Identity for the graded commutator of derivations, which in this case can be written as
$$
[D_1 ,[D_2 ,D_3 ]]=[[D_1 ,D_2 ],D_3 ]+(-1)^{|D_1 ||D_2 |}[D_2 ,[D_1 ,D_3 ]].
$$
From the graded curvature, other curvature tensors associated to the connection $\nnabla$ can be obtained, taking contractions with respect to other tensors. The case of contractions with respect to a graded metric (Ricci and scalar curvatures) is considered in \cite{MoS-97}. In the next section, we will see the possibilities that appear when contracting with respect to a graded symplectic form.

Note that $\nnabla_{D} $ can be extended as an operator to all of $\Omega^G (M)$; for example, if $\boldsymbol{\omega}$ is a graded $2-$form, then $\nnabla_{D} \boldsymbol{\omega}$ is again a $2-$form, given by its action on homogeneous derivations $D_1 ,D_2$ (of respective $\mathbb{Z}-$degrees $|D_1 |,|D_2 |$):
\begin{eqnarray*}
(\nnabla_{D} \boldsymbol{\omega} )(D_1 ,D_2 ) &=& D(\boldsymbol{\omega} (D_1 ,D_2 ))-(-1)^{|D|k}\boldsymbol{\omega} (\nnabla_{D}^0 D_1 ,D_2)\\
&-&(-1)^{(|D|+1)k}\boldsymbol{\omega} (\nnabla_{D}^1 D_1 ,D_2) \\
&-& (-1)^{|D|(k+|D_1 |)}\boldsymbol{\omega} (D_1 ,\nnabla_{D}^0 D_2)\\
&-&(-1)^{(|D|+1)(k+|D_1 |)}\boldsymbol{\omega} (D_1 ,\nnabla_{D}^1 D_2),
\end{eqnarray*}
where we have used the decomposition of the graded connection in its homogeneous components, $\nnabla=\nnabla^0 +\nnabla^1$.

We have given the definitions in an intrinsic manner, but when expressed in local super-coordinates $\{ x^i ,\theta_a \}$, the notions of connections and curvature coincide with those of \cite{Wit-92}.

As a further remark, let us note that, as in the non graded case, given a graded metric $G$ on $(M,\mathcal{A})$ (whose definition runs parallel to \ref{def-connection}), there exist a unique Levi-Civit\`a graded connection, which is symmetric and torsion-free (see \cite{MoS-96}, Theorem 4.2). From the explicit expression of this connection, it is apparent that it is even whenever $G$ is homogeneous (even or odd). However, also as in the non-graded case, there is not a canonical connection associated to a graded symplectic form.

With the notions developed in this section at hand, Definition \ref{def-fedosov} can be readily generalized.
\begin{definition}
Let $(M,\bigwedge \mathcal{E})$ be a split supermanifold, endowed with a graded symplectic form $\boldsymbol{\omega}$ and a graded symplectic connection $\nnabla$ (i.e., such that $\nnabla \boldsymbol{\omega} =0$). The triple $((M,\bigwedge \mathcal{E}),\nnabla,\boldsymbol{\omega} )$ is called a Fedosov supermanifold. 
\end{definition}

\section{Quantization and Fedosov supermanifolds}
\subsection{Poisson supermanifolds and deformation quantization}
Consider a supermanifold $(M,\bigwedge \mathcal{E})$ endowed with a super-Poisson bracket $[\![ \cdot ,\cdot ]\!]$. In attempting 
to develop the quantization deformation program in this context, following the ideas presented in Section \ref{sec-2} and taking 
into account the theory of Section \ref{sec-4}, one should seek for a series of bidifferential graded operators $\mathbf{C}_j$ 
acting on $\bigwedge \mathcal{E} [[\hbar ]]$, in a way such that the new product $\star$ on $\bigwedge \mathcal{E} [[\hbar ]]$ 
given by
$$
\alpha \star \beta =\sum^{\infty}_{j=0}\mathbf{C}_j (\alpha ,\beta)\hbar^j,
$$
has (among others) the properties that $\mathbf{C}_0 (\alpha ,\beta )=\alpha \wedge \beta$ (the exterior product on $\bigwedge \mathcal{E}$) and 
\begin{equation}\label{eq-c1}
\mathbf{C}_1 (\alpha ,\beta )-(-1)^{|\alpha ||\beta |}\mathbf{C}_1 (\beta ,\alpha )=i\hbar [\![ \alpha ,\beta ]\!]
\end{equation}
(recall Definition \ref{def-1}).\\
This observation is the starting point of the very interesting paper \cite{Bor-00}. In it, M. Bordemann considers only even Poisson 
brackets given by the Rothstein data $(\Lambda ,q,\nabla)$, where $\Lambda$ is the Poisson bivector associated to the Poisson 
bracket $\{ \cdot ,\cdot \}$ on $M$, $q$ is a metric on $E$ and $\nabla$ its Levi-Civit\`a connection. The next step in Fedosov's 
programme is to construct the Weyl bundle \eqref{bundle-v}, $\mathcal{V}=\Gamma W(M)\otimes \Omega (M)=\Gamma (W(M)\otimes T^* M)$, 
and here is where Bordemann departs from what it should be the analogue of Fedosov's construction in the super category, as instead 
of taking the full graded generalization of the product of bundles $W(M)\otimes T^* M$ he just introduces a set of anticommuting 
parameters by taking the product of $E$ with $W(M)\otimes T^* M$. According to the theory developed in Section \ref{sec-4}, 
the graded analogue of $\mathcal{V}$, let us call it $\mathcal{V}^G$, should be $\Gamma W^G (M)\otimes \Omega^G (M)$, where $W^G (M)$ 
should contain an appropriate generalization of the Fedosov's product \eqref{fed-prod}, and $\Omega^G (M)$ is the bundle of graded 
differential forms. Also, looking at equation \eqref{eq-delta}, it would seem enough to replace the tensor product $f\otimes g$, defined for functions $f,g\in \mathcal{C}^{\infty}(TM)$, by the graded tensor product in a way such that $f\otimes g (u,v)=(-1)^{|u||g|}f(u)g(v)$. These issues have been addressed in a recent paper by K. Bering \cite{Ber-09}, but there are still some others pending.\\
For instance, note that this time we should have ``superfunctions'' on the supertangent space $\mathcal{ST}(M,\mathcal{A})$ (not just graded $1-$forms). To fully generalize this setting, we should consider the ``super'' analogue of the space $\mathbb{R}$, which is the linear supermanifold $\mathbb{R}^{1|1}$ (at least, one should consider this space if does not want to lose information regarding the odd variables, see \cite{SV-09} for a justification for the use of $\mathbb{R}^{1|1}$ in the context of differential super-equations). Another question is the kind of supermanifold used. It seems that Bering employs DeWitt supermanifolds, although free use is made of the results by Rothstein. Due to the remarks done in Section \ref{sec-4} regarding tangent superspaces and supervector bundles and the fact that they appear without further ado (just mimicking the classical formulas with some added signs), it would be good if a clear geometric formulation of the results in \cite{Ber-09} could be done.\\
Anyway, it is surprising that the construction in \cite{Bor-00} recovers the original even Poisson bracket (apart from some constant factor) through the first order supercommutator \eqref{eq-c1}, in whose construction a lot of ``super'' structure has been discarded. There is some work in progress to elucidate these matters \cite{Vall-11}, but we can say that the deformation quantization approach in the graded manifold setting is not completely understood. Notice that we are leaving aside the study of deformation of \emph{odd} Poisson brackets, which is not done in either of the works aforementioned. However, their study should be quite natural: as shown in the Example \ref{ex-koszul}, from a Poisson (non graded) manifold, it is possible to construct in a very simple geometric way a supermanifold endowed with an odd Poisson bracket whose deformations may have an interesting geometric interpretation.

\subsection{The geometry of Fedosov supermanifolds}
The possibility of constructing a formal deformation quantization on a supermanifold through the Fedosov formalism, has led in a natural way to the study of the geometric properties of Fedosov supermanifolds. References on this topic are \cite{GL-04,LR-06,AL-09,MMV-09}.\\
Given a Fedosov supermanifold $((M,\bigwedge \mathcal{E}),\nnabla,\boldsymbol{\omega} )$, these properties are encoded in the symplectic curvature tensor
$$
R(D_1 ,D_2 ,D_3 ,D_4 )=\boldsymbol{\omega} (R^{\nnabla}(D_1 ,D_2 )D_3 ,D_4 ),
$$
and the main ones are just generalizations of the Bianchi's Identity and those of Proposition \ref{pro-propcurv}.
\begin{proposition}\label{pro-gradedcurv}
Let $((M,\bigwedge \mathcal{E}),\nnabla,\boldsymbol{\omega} )$ be a Fedosov supermanifold. Then, the graded symplectic curvature has the following properties, for all $D_1 ,D_2 ,D_3 ,D_4 \in \mathcal{X}^G (M)$,
\begin{enumerate}[(1)]
\item $R(D_1 ,D_2 ,D_3 ,D_4 )=-(-1)^{|D_1||D_2|}R(D_2 ,D_1 ,D_3 ,D_4 )$,
\item\label{item2} $R(D_1 ,D_2 ,D_3 ,D_4 )=(-1)^{|D_3||D_3|}R(D_1 ,D_1 ,D_4 ,D_3 )$,
\item\label{item3} $R^{\nnabla}(D_1 ,D_2 )D_3 +(-1)^{|D_1|(|D_2|+|D_3|)}R^{\nnabla}(D_2 ,D_3 )D_1 $\\
$+(-1)^{|D_3|(|D_1|+|D_2|)}R^{\nnabla}(D_3 ,D_1 )D_2 =0$.
\item\label{item4} $R(D_1 ,D_2 ,D_3 ,D_4 )+(-1)^{|D_4|(|D_1|+|D_2|+|D_3|)}R(D_4 ,D_1 ,D_2 ,D_3 )$\\
$+(-1)^{(|D_1|+|D_2|)(|D_3|+|D_4|)}R(D_3 ,D_4 ,D_1 ,D_2 )$\\
$+(-1)^{|D_1|(|D_2|+|D_3|+|D_4|)}R(D_2 ,D_3 ,D_4 ,D_1 )=0$.
\end{enumerate}
\end{proposition}
\noindent Note that, as a consequence of \eqref{item3}, we also have \\
\begin{eqnarray*}
R(D,D_1 ,D_2 ,D_3 ) &+& (-1)^{|D_1|(|D_2|+|D_3|)}R(D,D_2 ,D_3 ,D_1 ) \\
&+&(- 1)^{|D_3|(|D_1|+|D_2|)}R(D,D_3 ,D_1 ,D_2)=0.
\end{eqnarray*}

Indeed, this result is used (permuting arguments and adding up the resulting equations) to prove \eqref{item4}. The first item is immediate from the definition, and the items \eqref{item2} and \eqref{item3} are proved exactly in the same way as in the non graded case.

Now, generalizing \ref{teo-contractions}, we define the graded symplectic Ricci tensor.
\begin{definition}
Let $((M,\bigwedge \mathcal{E}),\nnabla,\boldsymbol{\omega} )$ be a Fedosov supermanifold. Its graded symplectic Ricci tensor $\boldsymbol{\rho}$ is defined, for any $D_1 ,D_2 \in \mathcal{X}^G (M)$, as
$$
\boldsymbol{\rho}(D_1 ,D_2 )=\mathcal{S}\mathrm{Tr}_{\boldsymbol{\omega}}(D\mapsto R^{\nnabla}(D,D_1 )D_2),
$$
where $\mathcal{S}\mathrm{Tr}_{\boldsymbol{\omega}}$ is the supertrace with respect to $\boldsymbol{\omega}$ considered as a bilinear superform \cite{Lei-80}.
\end{definition}
By using the properties stated in \ref{pro-gradedcurv}, it is easy to see that, as in the classical case, there is basically only one possible contraction to define the graded symplectic Ricci tensor, the other possibilities give a vanishing tensor or a multiple of $\boldsymbol{\rho}$, depending on the parity of $\boldsymbol{\omega}$. What is really interesting is the phenomenon that appears when taking a further contraction of the Ricci tensor with the graded symplectic form $\boldsymbol{\omega}$. In the classical case, Corollary \ref{cor-szero}, this construction gives zero by the opposite symmetry properties of $\boldsymbol{\omega}$ and $\boldsymbol{\rho}$, and this is the case again for \emph{even} symplectic forms, as in this case $\boldsymbol{\rho}$ is graded symmetric:
$$
\boldsymbol{\rho}(D_1 ,D_2 )=(-1)^{|D_1||D_2|}\boldsymbol{\rho}(D_2 ,D_1 ).
$$
However, the corresponding $\boldsymbol{\rho}$ for an \emph{odd} symplectic form $\boldsymbol{\omega}$, does not have definite symmetry properties (it is neither graded symmetric or skew-symmetric), so when defining the odd scalar curvature $S$ as the contraction of $\boldsymbol{\rho}$ with an odd symplectic form $\boldsymbol{\omega}$, it results that $S\neq 0$. In \cite{BB-08,BB-08b}, I. Batalin and K. Bering give an interesting physical interpretation of this scalar curvature, but there are some open questions regarding its geometrical origin. To discuss one of them, note that odd symplectic forms of second order depth on a split supermanifold $(M,\bigwedge \mathcal{E})$ are particularly simple, and can be described by a tensor field on the base manifold $M$. This is due to the following result \cite{Mon-92}.
\begin{theorem}
Let $\nabla$ be a linear connection on the vector bundle $E$ (of rank $n$) and let $\boldsymbol{\omega}$ be a closed graded $2-$form on $(M,\bigwedge \mathcal{E})$. Then, $\boldsymbol{\omega}$ is uniquely determined by the following tensor fields:
\begin{enumerate}[(i)]
\item A closed ordinary $2-$form $\tilde{\omega}\in \Omega (M)$,
\item $\bar{K}=\sum_{k=1}^n K_k \in \sum_{k=1}^n \Gamma (T^* M\otimes \Lambda^k E)$,
\item $\bar{L}=\sum_{k=2}^n L_k \in \sum_{k=1}^n \Gamma^s (E\otimes \Gamma^{k-1}E$,
\end{enumerate}
where $\tilde{\omega}, K_1$ and $L_2$ are independent of the chosen connection $\nabla$.
\end{theorem}
Note that, for second order depth forms, the sums run up only to $k=2$. A corollary to this result is the following.
\begin{corollary}
If $\boldsymbol{\omega}$ is an odd graded symplectic form, then $\tilde{\omega}=0=L_2$, and $K_1$ defines a bundle isomorphism between $TM$ and $E$.
\end{corollary}
Thus, for example, starting only from a usual symplectic form $\tilde{\omega}\in \Omega^2 (M)$, one could construct an odd symplectic form $\boldsymbol{\omega}$ on the geometric supermanifold $(M,\Omega (M)$ (taking $K_1:TM\to T^* M$ as the inverse Poisson map associated to $\tilde{\omega}$). The results in \cite{MMV-09} provide a way to construct also a symplectic connection $\nnabla$ from $\tilde{\omega}$, and the resulting odd symplectic scalar curvature is solely determined in terms of geometric objects defined on $M$. These ideas may be useful to study the geometric meaning of the odd symplectic scalar curvature and its relationship to the physical interpretation presented in \cite{BB-08b}.

\subsection*{Acknowledgements}
Work partially supported by a CONACyT grant code CB-J2-78791 (M\'exico).

\end{document}